\documentclass[a4paper,11pt,reqno,twoside]{article}
\usepackage{amssymb}
\usepackage[frenchb,english]{babel}
\usepackage{amscd}
\usepackage{amsmath,amsthm,amsfonts,amssymb,graphicx}
\usepackage[CJKbookmarks,colorlinks,linkcolor=blue,anchorcolor=blue,citecolor=green]{hyperref}
\usepackage{mathrsfs}
\usepackage{fancyhdr}
\usepackage{subfigure}
\usepackage{bm}
\usepackage{multicol}
\usepackage{picins}
\usepackage{abstract}

\makeatother
\begin{document}
\theoremstyle{plain}
\newtheorem{Definition}{\hspace{0em}{\bf{Definition}}}[section]
\newtheorem{Proposition}{\hspace{0em}{\bf{Proposition}}}[section]
\newtheorem{Property}{\hspace{0em}{\bf{Property}}}[section]
\newtheorem{Theorem}{\hspace{0em}\bf{Theorem}}[section]
\newtheorem{Lemma}[Theorem]{\hspace{0em}\bf{Lemma}}
\newtheorem{Corollary}[Theorem]{\hspace{0em}{\bf{Corollary}}}
\newtheorem{Axiom}{\hspace{0em}{\bf{Axiom}}}[section]
\newtheorem{Exercise}{\hspace{0em}{\bf{Exercise}}}[section]
\newtheorem{Question}{\hspace{0em}{\bf{Question}}}
\newtheorem{Example}{\hspace{0em}{\bf{Example}}}
\newtheorem{Notation}{\hspace{0em}{\bf{Notation}}}
\newtheorem{Remark}{\hspace{0em}{\bf{Remark}}}[section]
%%%%%%%%%%%%%%%%%%%%%%%%%%%%%%%%%%%%%%%%%%%%%%%%%%%%%%%%%%%%%%%%

%%%%%%%%%%%%%%%%%%%%%%%%%%%%%%%%%%%%%%%%%%%%%%%%%%%%%%%%%%%%%%%%
\setlength{\oddsidemargin}{ 1cm}  % 3.17cm - 1 inch
\setlength{\evensidemargin}{\oddsidemargin}
\setlength{\textwidth}{13.50cm}
\vspace{-.8cm}

\noindent  {\LARGE On canonical metrics on Cartan-Hartogs domains}\\\\
\noindent\text{Zhiming Feng  }\\
\noindent\small {School of Mathematical and Information Sciences, Leshan Normal University, Leshan, Sichuan 614000, P.R. China } \\
\noindent\text{Email: fengzm2008@163.com}

\vskip 5pt
\noindent\text{Zhenhan Tu$^{*}$ }\\
\noindent\small {School of Mathematics and Statistics, Wuhan
University, Wuhan, Hubei 430072, P.R. China} \\
\noindent\text{Email: zhhtu.math@whu.edu.cn
}
\renewcommand{\thefootnote}{{}}
\footnote{\hskip -16pt {$^{*}$Corresponding author. \\}}
\\

\normalsize \noindent\textbf{Abstract}\quad {The Cartan-Hartogs
domains are defined as a class of Hartogs type domains over
irreducible bounded symmetric domains. The purpose of this paper is
twofold. Firstly, for a Cartan-Hartogs domain
$\Omega^{B^{d_0}}(\mu)$ endowed with the canonical metric $g(\mu),$
we obtain an explicit formula for the Bergman kernel of the weighted
Hilbert space $\mathcal{H}_{\alpha}$ of square integrable
holomorphic functions on $(\Omega^{B^{d_0}}(\mu), g(\mu))$ with the
weight $\exp\{-\alpha \varphi\}$ (where $\varphi$ is a globally
defined K\"{a}hler potential for $g(\mu)$) for $\alpha>0$, and,
furthermore, we give an explicit expression of the Rawnsley's
$\varepsilon$-function expansion for $(\Omega^{B^{d_0}}(\mu),
g(\mu)).$ Secondly, using the explicit expression of the Rawnsley's
$\varepsilon$-function expansion, we show that the coefficient $a_2$
of the Rawnsley's $\varepsilon$-function expansion for the
Cartan-Hartogs domain $(\Omega^{B^{d_0}}(\mu), g(\mu))$ is constant
on $\Omega^{B^{d_0}}(\mu)$ if and only if $(\Omega^{B^{d_0}}(\mu),
g(\mu))$ is biholomorphically isometric to the complex hyperbolic
space. So we give an affirmative answer to a conjecture raised by M.
Zedda.
\smallskip\\\
\textperiodcentered \; K\"{a}hler metrics
\smallskip\\\
  \textperiodcentered \, 32M15  \textperiodcentered \, 32Q15
%%%%%%%%%%%%%%%%%%%%%%%%%%%%%%%%%%%%%%%%%%%%%%%%%%%%%%%%%%%%%%%%

%%%%%%%%%%%%%%%%%%%%%%%%%%%%%%%%%%%%%%%%%%%%%%%%%%%%%%%%%%%%%%%%
\setlength{\oddsidemargin}{-.5cm}  % 3.17cm - 1 inch
\setlength{\evensidemargin}{\oddsidemargin}
\pagenumbering{arabic}
\renewcommand{\theequation}
{\arabic{section}.\arabic{equation}}
%%%%%%%%%%%%%%%%%%%%%%%%%%%%%%%%%%%%%%%%%%%%%%%%%%%%%%%%%%%%%%%%

 \setcounter{equation}{0}
\section{{Introduction}}
The expansion of the Bergman kernel has received a lot of attention
recently, due to the influential work of Donaldson, see e.g.
\cite{Donaldson}, about the existence and uniqueness of constant
scalar curvature K\"{a}hler metrics (cscK metrics). Donaldson used
the asymptotics of the Bergman kernel proved by Catlin \cite{Cat}
and Zelditch \cite{Zeld} and the calculation of Lu \cite{Lu} of the
first coefficient in the expansion to give conditions for the
existence of cscK metrics. This work inspired many papers on the
subject since then. For the reference of the expansion of the
Bergman kernel, see also Engli\v{s} \cite{E2}, Loi \cite{Lo},
Ma-Marinescu \cite{MM07, MM08, MM12}, Xu \cite{X}
 and references therein.

Assume that $D$ is a bounded domain in $\mathbb{C}^n$ and $\varphi$
is a strictly plurisubharmonic function on $D$. Let $g$ be a
K\"{a}hler  metric  on $D$ associated to the K\"{a}hler form
$\omega=\frac{\sqrt{-1}}{2\pi}\partial\overline{\partial}\varphi$.
For $\alpha>0$, let $\mathcal{H}_{\alpha}$ be the weighted Hilbert
space of square integrable holomorphic functions on $(D, g)$ with
the weight $\exp\{-\alpha \varphi\}$, that is,
$$\mathcal{H}_{\alpha}:=\left\{ f\in \textmd{Hol}(D)\;\left |\; \int_{D}\right.|f|^2\exp\{-\alpha \varphi\}\frac{\omega^n}{n!}<+\infty\right\},$$
where $\textmd{Hol}(D)$ denotes the space of holomorphic functions
on $D$. Let $K_{\alpha}$ be the Bergman kernel (namely, the
reproducing kernel) of $\mathcal{H}_{\alpha}$ if
$\mathcal{H}_{\alpha}\neq \{0\}$.  The Rawnsley's
$\varepsilon$-function on $D$ (see Cahen-Gutt-Rawnsley \cite{CGR}
and
 Rawnsley \cite{R}) associated to the metric $g$  is defined by
\begin{equation}\label{eq1.4}
 \varepsilon_{\alpha}(z):=\exp\{-\alpha \varphi(z)\}K_{\alpha}(z,\overline{z}),\;\; z\in D.
\end{equation}
Note the Rawnsley's $\varepsilon$-function depends only on the
metric $g$ and not on the choice of the K\"{a}hler potential
$\varphi$ (which is defined up to an addition with the real part of
a holomorphic function on $D$). If the function
$\varepsilon_{\alpha}(z)$ $(z\in D)$ is a positive constant for
$\alpha=1$, the metric  $g$ on $D$ is called to be balanced.

The asymptotics of $ \varepsilon_{\alpha}$ was expressed in terms of
the parameter $\alpha$ for compact manifolds by Catlin \cite{Cat}
and Zelditch \cite{Zeld} (for $\alpha\in \mathbb{N}$) and for
non-compact manifolds by Ma-Marinescu \cite{MM07,MM08}. In some
particular case it was also proved by Engli\v{s} \cite{E1,E2}.

The Cartan-Hartogs domains are defined as a class of Hartogs type
domains over irreducible bounded symmetric domains. Let $\Omega$ be
an irreducible bounded symmetric domain in $\mathbb{C}^d$ of genus
$p$. The generic norm of $\Omega$ is defined by $ N(z, \bar \xi) :=
(V(\Omega ) K(z, \bar \xi))^{-1/p },$ where $V(\Omega )$ is the
total volume of $\Omega$ with respect to the Euclidean measure of
$\mathbb{C}^d$ and $K(z, \overline{\xi})$ is its Bergman kernel. For
an irreducible bounded symmetric domain $\Omega$ in $\mathbb{C}^d$,
a positive real number $\mu$ and a positive integer number $d_0$,
the Cartan-Hartogs domain $\Omega^{B^{d_0}}(\mu)$ is defined by
\begin{equation}\label{eq1.1}
  \Omega^{B^{d_0}}(\mu):=\left.\left\{(z,w)\in  \Omega\times \mathbb{C}^{d_0}\subset \mathbb{C}^{d}\times \mathbb{C}^{d_0}
  ~\right|~ \|w\|^2<N(z,\overline{z})^{\mu} \right\},
\end{equation}
where $\|\cdot\|$ is the standard Hermitian norm in
$\mathbb{C}^{d_0}$.

Let $\mathcal{M}_{m,n}$ be the set of all $m\times n$ matrices
$z=(z_{ij})$ with complex entries. Let ${\overline z}$ be the
complex conjugate of the matrix $z$ and let ${z}^t$ be the transpose
of the matrix $z$. $I$ denotes the identity matrix. If a square
matrix $z$ is positive definite, then we write $z>0$. For each
bounded classical  symmetric domain $\Omega$ (refer to Hua
\cite{Hua}), we list the genus $p(\Omega)$, the generic norm
$N_\Omega(z,\overline{z})$ of $\Omega$ and corresponding
Cartan-Hartogs domain $\Omega^{B^{d_0}}(\mu)$ \cite{Yin} according
to its type as following.

 $(i)$ If $\Omega=\Omega_I(m,n):=\{z\in
\mathcal{M}_{m,n}: I-z{\overline z}^t>0 \}$ ($1\leq m\leq n$) (the
classical domains of type $I$), then $p(\Omega)=m+n,$
$N_\Omega(z,\overline{z})=\det(I-z{\overline z}^t),$ and
$$\Omega^{B^{{d_0}}}(\mu)=\left\{(z,w)\in \Omega_I(m,n)\times
\mathbb{C}^{{d_0}}\subset \mathbb{C}^{mn}\times \mathbb{C}^{{d_0}}:
 \|w\|^{{2}}<  (\det(I-z{\overline z}^t) )^{\mu}  \right\}.$$
Specially, when $\Omega:=B^n$ is the unit ball in $ \mathbb{C}^{n}$,
then we have
$$\Omega^{B^{{d_0}}}(\mu)
   =\left\{(z,w)\in \Omega\times \mathbb{C}^{{d_0}}:
  \|z\|^2+ \|w\|^{\frac{2}{\mu}} <1 \right\}.$$
It is a natural generalization of Thullen domains.

$(ii)$ If $\Omega= \Omega_{II}(n):=\{z\in \mathcal{M}_{n,n}: z^t=-z,
I-z{\overline z}^t>0 \}$ ($n\geq 4$) (the classical domains of type $II$), then
$p(\Omega)=2(n-1),$
$N_\Omega(z,\overline{z})=(\det(I-z{\overline z}^t))^{1/2} ,$ and
$$\Omega^{B^{{d_0}}}(\mu)=\left\{(z,w)\in  \Omega_{II}(n)\times
\mathbb{C}^{{d_0}}\subset \mathbb{C}^{n(n-1)/2}\times
\mathbb{C}^{{d_0}}:
 \|w\|^{{2}}<  (\det(I-z{\overline z}^t) )^{\mu/2}  \right\}.$$

$(iii)$ If $\Omega= \Omega_{III}(n):=\{z\in \mathcal{M}_{n,n}: z^t=z,
I-z{\overline z}^t>0 \}$ ($n\geq 2$) (the classical domains of type $III$), then
$p(\Omega)=n+1,$   $N_\Omega(z,\overline{z})=\det(I-z{\overline
z}^t) ,$ and
$$\Omega^{B^{{d_0}}}(\mu)=\left\{(z,w)\in \Omega_{III}(n)\times
\mathbb{C}^{{d_0}}\subset \mathbb{C}^{n(n+1)/2}\times
\mathbb{C}^{{d_0}}:
 \|w\|^{{2}}<  (\det(I-z{\overline z}^t) )^{\mu}  \right\}.$$

$(iv)$ If $\Omega=\Omega_{IV}(n):=\{z\in \mathbb{C}^{n}: 1-2z{\overline
z}^t+|zz^t|^2>0,  z{\overline z}^t<1\} $ ($n\geq 5$) (the classical domains of
type $IV$), then
 $p(\Omega)=n,$
$N_\Omega(z,\overline{z})=1-2z{\overline z}^t +|zz^t|^2,$ and
$$\Omega^{B^{{d_0}}}(\mu)=\left\{(z,w)\in  \Omega_{IV}(n)\times
\mathbb{C}^{{d_0}}\subset \mathbb{C}^{n}\times \mathbb{C}^{{d_0}}:
 \|w\|^{{2}}<  (1-2z{\overline z}^t +|zz^t|^2)^{\mu}  \right\}.$$

For the Cartan-Hartogs domain $\Omega^{B^{d_0}}(\mu)$, define
\begin{equation}\label{eq1.2}
 \Phi(z,w):=-\log(N(z,\overline{z})^{\mu}-\|w\|^2).
\end{equation}
The K\"{a}hler form $\omega(\mu)$ on  $\Omega^{B^{d_0}}(\mu)$
is defined by
\begin{equation}\label{eq1.3}
  \omega(\mu):=\frac{\sqrt{-1}}{2\pi}\partial
\overline{\partial}\Phi.
\end{equation}
 The K\"{a}hler metric $g(\mu)$ on  $\Omega^{B^{d_0}}(\mu)$
associated to  $\omega(\mu)$ is given by
$ds^2=\sum_{i,j=1}^{n}\frac{\partial^2\Phi}{\partial z_i \partial
\overline{z_j}}dz_i\otimes d\overline{z_j},$ where $n=d+d_0$,
$z=(z_1,z_2,\cdots,z_d)$, $w=(z_{d+1},z_{d+2},\cdots,z_n)$. With the
exception of the complex hyperbolic space which is obviously
homogeneous, each Cartan-Hartogs domain $(\Omega^{B^{d_0}}(\mu),
g(\mu))$ is a noncompact, nonhomogeneous, complete K\"{a}hler
manifold (see Yin-Wang \cite{YW}). Further, for some particular
value $\mu_0$ of $\mu$, $g(\mu_0)$ is a K\"{a}hler-Einstein metric.
For the general reference of the Cartan-Hartogs domains in this
paper, see Loi-Zedda \cite{LZ}, Wang-Yin-Zhang-Roos \cite{RWYZ}, Yin
\cite{Yin}, Yin-Wang \cite{YW}, Zedda \cite{Zed} and references
therein.

In this paper, we study the asymptotics of the Rawnsley's
$\varepsilon$-function on the Cartan-Hartogs domain with the
canonical metric and draw some geometric consequences. For a
Cartan-Hartogs domain $(\Omega^{B^{d_0}}(\mu), g(\mu)),$ we have
(see Theorem 3.1 in this paper) that the Rawnsley's
$\varepsilon$-function admits the expansion:
\begin{equation}\label{eq1.7}
  \varepsilon_{\alpha}(z,w)=\sum_{j=0}^{d+d_0}a_j(z,w)\alpha^{d+d_0-j}, \;\; (z,w)\in
  \Omega^{B^{d_0}}(\mu).
\end{equation}

By Th. 1.1 of Lu \cite{Lu},  Th. 4.1.2 and Th. 6.1.1 of Ma-Marinescu
\cite{MM07}, Th. 3.11 of Ma-Marinescu \cite{MM08} and Th. 0.1 of
Ma-Marinescu \cite{MM12}, see also Th. 3.3 of Xu \cite{X}, we
have
\begin{equation}\label{eq1.6}
\left\{  \begin{array}{ll}
    a_0 & =1, \\
    a_1 & = \frac{1}{2}k_g, \\
    a_2 & =\frac{1}{3}\triangle k_g+\frac{1}{24}|R|^2-\frac{1}{6}|Ric|^2+\frac{1}{8}k_g^2,
  \end{array}\right.
\end{equation}
where $k_g$, $\triangle$, $R$ and $Ric$  denote  the scalar
curvature, the Laplace, the curvature tensor and the Ricci curvature
associated to the metric $g(\mu)$, respectively.

Let $ {B}^{d}:=\left\{z=(z_1,z_2,\cdots,z_d)\in \mathbb{C}^d\;\left
|\; \|z\|^2\right.=\sum_{k=1}^d|z_k|^2<1\right\}$ and  let the
metric $g_{hyp}$ on  $ {B}^{d}$  be given by $
ds^2=-\sum_{i,j=1}^{d}\frac{\partial^2\log(1-\|z\|^2)}{\partial z_i
\partial \overline{z_j}}dz_i\otimes d\overline{z_j}.$ Then we call $(B^d, g_{hyp})$
the complex hyperbolic space. Note that here
$\mathcal{H}_{\alpha}\neq \{0\}$ iff $\alpha > d$ and that $\alpha
g_{hyp}\;(\alpha>0)$ is a balanced metric on $B^d$ iff $\alpha>d$.

Loi and Zedda \cite{LZ} studied balanced metrics on the
Cartan-Hartogs domain and proved the following result for $d_0=1$:

\begin{Theorem}\textup{(Loi-Zedda \cite{LZ} for $d_0=1$)}\label{Th:a.1.3}{
Let $\Omega$ be an irreducible bounded symmetric domain of dimension
$d$ and genus $p$.  Then the metric  $\alpha g(\mu)$ on
$\Omega^{B^{d_0}}(\mu)$ is balanced if and only if
$\alpha>\max\{d+d_0,\frac{p-1}{\mu}\}$ and $(\Omega^{B^{d_0}}(\mu),
g(\mu))$ is holomorphically isometric to the complex hyperbolic
space $({B}^{d+d_0}, g_{hyp})$, namely, $\Omega=B^d$ and  $\mu=1$ .
 }\end{Theorem}

By calculating the scalar curvature $k_g$, the Laplace $\Delta k_g$
of $k_g$, the norm $|R|^2$ of the curvature tensor $R$ and the norm
$|Ric|^2$ of the Ricci curvature $Ric$ of a Cartan-Hartogs domain
$(\Omega^{B^{d_0}}(\mu), g(\mu))$, Zedda \cite{Zed} has proved the
following theorem for $d_0=1$:

\begin{Theorem}\textup{(Zedda \cite{Zed} for $d_0=1$)}\label{Th:a.1.1}{
Let $(\Omega^{B^{d_0}}(\mu), g(\mu))$ be a Cartan-Hartogs domain. If
the coefficient $a_2$ of the Rawnsley's $\varepsilon$-function
expansion is a constant on $\Omega^{B^{d_0}}(\mu)$, then
$(\Omega^{B^{d_0}}(\mu), g(\mu))$ is K\"{a}hler-Einstein.
 }\end{Theorem}
Further, Zedda \cite{Zed} conjectured that the coefficient $a_2$ of
the expansion of the Rawnsley's $\varepsilon$-function associated to
$g(\mu)$ is constant if and only if $(\Omega^{B^{d_0}}(\mu),
g(\mu))$ is biholomorphically isometric to the complex hyperbolic
space. Obviously, the conjecture implies Theorem 1.2.  In this
paper, for any positive integer $d_0$, by giving an explicit
expression of the reproducing kernel $K_{\alpha}$ of
$\mathcal{H}_{\alpha}$ and the Rawnsley's $\varepsilon$-function for
$(\Omega^{B^{d_0}}(\mu), g(\mu))$, we prove that the Zedda's
 conjecture is affirmative, namely, we prove the
following conclusion:

\begin{Theorem}\label{Th:a.1.2}{Let $(\Omega^{B^{d_0}}(\mu), g(\mu))$ be a Cartan-Hartogs domain.
 Then the  coefficient $a_2$ of the  Rawnsley's
$\varepsilon$-function expansion is a constant on
$\Omega^{B^{d_0}}(\mu)$ if and only if $(\Omega^{B^{d_0}}(\mu),
g(\mu))$ is biholomorphically isometric to the complex hyperbolic
space $({B}^{d+d_0}, g_{hyp})$.
 }\end{Theorem}

Remark that Theorem \ref{Th:a.1.2} imediately implies Theorem
\ref{Th:a.1.3}. In fact, for $\alpha>\max\{d+d_0,\frac{p-1}{\mu}\}$,
let $\omega_{\alpha}:=\frac{\sqrt{-1}}{2\pi}\partial
\overline{\partial}(-\alpha\log(N(z,\overline{z})^{\mu}-\|w\|^2))
=\frac{\sqrt{-1}}{2\pi}\partial \overline{\partial}(\alpha \Phi)$ (so $
\omega_{\alpha}= \alpha \omega{(\mu)}$),
$$\mathcal{H}_{\omega_{\alpha}}:=\left\{ f\in
\textmd{Hol}(\Omega^{B^{d_0}}(\mu))\;\left|\;
\int_{\Omega^{B^{d_0}}(\mu)}\right.|f|^2\exp\{-\alpha
\Phi\}\frac{\omega_{\alpha}^{d+d_0}}{(d+d_0)!}<+\infty\right\},$$
and
$$\mathcal{H}_{{\alpha}}:=\left\{ f\in
\textmd{Hol}(\Omega^{B^{d_0}}(\mu))\;\left|\;
\int_{\Omega^{B^{d_0}}(\mu)}\right.|f|^2\exp\{-\alpha
\Phi\}\frac{\omega(\mu)^{d+d_0}}{(d+d_0)!}<+\infty\right\}.$$ It is
easy to see that
$K_{\omega_{\alpha}}=\frac{1}{\alpha^{d+d_0}}K_{\alpha}$, where
$K_{\omega_{\alpha}}$ and $K_{\alpha}$ are the Bergman kernels of
$\mathcal{H}_{\omega_{\alpha}}$ and $\mathcal{H}_{{\alpha}}$,
respectively. So, we have
$\exp\{-\alpha\Phi\}K_{\omega_{\alpha}}=\frac{1}{\alpha^{d+d_0}}\exp\{-\alpha\Phi\}K_{\alpha}$.
By the definition, $\alpha g(\mu)$ on $\Omega^{B^{d_0}}(\mu)$ is
balanced iff $\exp\{-\alpha\Phi\}K_{\omega_{\alpha}}$ is a positive
constant on $\Omega^{B^{d_0}}(\mu)$. This indicate that $\alpha
g(\mu)$ on $\Omega^{B^{d_0}}(\mu)$ is balanced if and only if
$\varepsilon_{\alpha}=\exp\{-\alpha\Phi\}K_{\alpha}$ is a positive
constant on $\Omega^{B^{d_0}}(\mu)$. Note that $\alpha g_{hyp}$ is a
balanced metric on $B^d$ iff $\alpha>d$, and if the metric $\alpha
g(\mu)$ on $\Omega^{B^{d_0}}(\mu)$ is balanced, then we have
$\alpha>\max\{d+d_0,\frac{p-1}{\mu}\}$ (see Lemma 9 in \cite{LZ}).
Thus, by Theorem \ref{Th:a.1.2}, we obtain that the metric $\alpha
g(\mu)$ on $\Omega^{B^{d_0}}(\mu)$ is balanced if and only if
$\alpha>\max\{d+d_0,\frac{p-1}{\mu}\}$ and $(\Omega^{B^{d_0}}(\mu),
g(\mu))$  is holomorphically isometric to the complex hyperbolic
space $({B}^{d+d_0}, g_{hyp})$. The proof is complete.

On the other hand, combining the formulas \eqref{eq2.38} and
\eqref{eq2.61} in this paper, we have that $a_1$ is constant if and
only if $\mu=\frac{p}{d+1}$. Further, by \eqref{eq1.6}, the scalar
curvature of a Cartan-Hartogs domain
$(\Omega^{B^{d_0}}(\mu),g(\mu))$ is constant iff $a_1$ is constant.
Thus we get the following theorem:

\begin{Theorem}\textup{(Zedda \cite{Zed} for $d_0=1$)}\label{Th:a.1.4}{
Let $\Omega$ be an irreducible bounded symmetric domain of dimension
$d$ and genus $p$. Then the scalar curvature of a Cartan-Hartogs
domain $(\Omega^{B^{d_0}}(\mu),g(\mu))$ is constant if and only if
$g(\mu)$ is K\"{a}hler-Einstein, namely, $\mu=\frac{p}{d+1}$.
 }\end{Theorem}

The paper is organized as follows. In Section 2, we obtain an
explicit formula for the Bergman kernel $K_{\alpha}$ of
$\mathcal{H}_{\alpha}$ for the Cartan-Hartogs domain
$(\Omega^{B^{d_0}}(\mu),g(\mu))$ in terms of ranks, Hua polynomials
and generic norms of $\Omega$ and $B^{d_0}$ (Theorem 2.3). In
Section 3, using results in Section 2, we give the explicit
expansion of the Rawnsley's $\varepsilon$-function and obtain the
expression of its coefficients $a_1, a_2$ for the Cartan-Hartogs
domain associated to $g(\mu)$ (Corollary 3.2). Finally, in Section
4, the conclusion is achieved by using the classification of bounded
symmetric domains (it follows that $a_2$ is constant if and only if
the rank $r = 1$ and $\mu= 1$).

\setcounter{equation}{0}
\section{The reproducing kernel of $\mathcal{H}_{\alpha}$ for
$\Omega^{B^{d_0}}(\mu)$  with the canonical metric $g(\mu)$}

Let $\Omega$ be an irreducible bounded symmetric domain  in
$\mathbb{C}^d$ in its Harish-Chandra realization. Thus $\Omega$ is
the open unit ball of a Banach space which admits the structure of a
$JB^{\ast}$-triple. We denote  $r, a, b, d, p $ and
$N(z,\overline{w})$ by the rank, the characteristic multiplicities,
the dimension,  the genus, and the generic norm of $\Omega$,
respectively. Thus
\begin{equation}\label{1.1}
    d=\frac{r(r-1)}{2}  a+rb+r,\quad   p=(r-1)a+b+2.
\end{equation}
For any $s>-1$, the value of the Hua integral
$\int_{\Omega}N(z,\overline{z})^s dm(z)$ is given by
\begin{equation}\label{1.2}
\int_{\Omega}N(z,\overline{z})^s
dm(z)=\frac{\chi(0)}{\chi(s)}\int_{\Omega}dm(z),
\end{equation}
where $dm(z)$ denotes the Euclidean measure on $\mathbb{C}^d$, $\chi$
is the Hua polynomial
\begin{equation}\label{1.3}
   \chi(s):=\prod_{j=1}^r\left(s+1+(j-1)\frac{a}{2}\right
   )_{1+b+(r-j)a},
\end{equation}
in which, for a  non-negative integer $m$, $(s)_m$ denotes the
raising factorial
$$ {(s)_m:=\frac{\Gamma(s+m)}{\Gamma(s)}=s(s+1)\cdots (s+m-1)}.$$

Let $\mathcal{G}$ stand for the identity connected component of the
group of biholomorphic self-maps of $\Omega$, and $\mathcal{K}$ {for
the stabilizer} of the origin in $\mathcal{G}$. Under the action
$f\mapsto f\circ k (k\in \mathcal{K})$ of $\mathcal{K}$, the space
$\mathcal{P}$ of holomorphic polynomials on $\mathbb{C}^d$ admits
the Peter-Weyl decomposition
$$\mathcal{P}=\bigoplus_{\lambda}\mathcal{P}_{\lambda},$$
 {where the summation is taken over all partitions}
$\lambda$, i.e., $r-$tuples $(\lambda_1, \lambda_2, \cdots,
\lambda_r)$ of nonnegative integers such that $\lambda_1\geq
\lambda_2 \geq \cdots \geq \lambda_r \geq 0$,  {and the spaces}
$\mathcal{P}_{\lambda}$ are $\mathcal{K}$-invariant and irreducible.
For each $\lambda$, $\mathcal{P}_{\lambda} \subset
\mathcal{P}_{|\lambda|}$, where $|\lambda|$ denotes the weight of
partition $\lambda$, i.e., $|\lambda|:=\sum_{j=1}^r \lambda_j$,  and
$\mathcal{P}_{|\lambda|}$ is the space of homogeneous holomorphic
polynomials of degree $|\lambda|$.

Let
\begin{equation}\label{1.4}
    {\langle}f,g {\rangle}_{\mathcal{F}}:=\int_{\mathbb{C}^d}f(z)\overline{g(z)} d\rho_{\mathcal{F}}(z)
\end{equation}
be the Fock-Fischer inner product on the space $\mathcal{P}$ of
holomorphic polynomials on $\mathbb{C}^d$, where
\begin{equation}\label{1.5}
    d\rho_{\mathcal{F}}(z):=\frac{1}{\pi^d}e^{-\|z\|^2}dm(z).
\end{equation}
 For every partition $\lambda$, let $K_{\lambda}(z_1,\overline{z_2})$ be
the Bergman kernel of $\mathcal{P}_{\lambda}$ with respect to
\eqref{1.4}. The weighted Bergman kernel of the weighted Hilbert
space $A^2(\mathbb{C}^d,\rho_{\mathcal{F}})$ of square-integrable
holomorphic functions on $\mathbb{C}^d$ with the measure
$d\rho_{\mathcal{F}}$ is
\begin{equation}\label{1.6}
    K(z_1,\overline{z_2}):=\sum_{\lambda}K_{\lambda}(z_1,\overline{z_2}).
\end{equation}

The kernels $K_{\lambda}(z_1,\overline{z_2})$ are related to the
generic norm  $N(z_1,\overline{z_2})$ by the Faraut-{Kor\'{a}nyi}
formula
\begin{equation}\label{1.7}
N(z_1,\overline{z_2})^{-s}=\sum_{\lambda}(s)_{\lambda}K_{\lambda}(z_1,\overline{z_2}),
\end{equation}
the series converges  {uniformly} on compact subsets of
$\Omega\times\Omega$, $s\in \mathbb{C}$,  where $(s)_{\lambda}$
denote the generalized Pochhammer symbol
\begin{equation}\label{1.8}
   (s)_{\lambda}:=\prod_{j=1}^r\left(s-\frac{j-1}{2}a\right )_{\lambda_j}.
\end{equation}
For the proofs of above facts and additional details, we refer,
e.g.,  to \cite{FK}, \cite{FKKLR} and \cite{YLR}.

\begin{Lemma}\label{Le:2.1}{
Let $\Omega$ be an irreducible bounded symmetric domain  in
$\mathbb{C}^d$ in its Harish-Chandra realization with the generic
norm $N$ and the genus $p$. For $z_0\in \Omega$, let $\phi$ be an
automorphism of  $\Omega$ such that $\phi(z_0)=0$. By \cite{RWYZ},
the function
\begin{equation}\label{eq2.1}
   \psi(z):=\frac{N(z_0,\overline{z_0})^{\frac{\mu}{2}}}{N(z,\overline{z_0})^{\mu}}
\end{equation}
satisfies
\begin{equation}\label{eq2.2}
    |\psi(z)|^2=\left(\frac{N(\phi(z),\overline{\phi(z)})}{N(z,\overline{z})}\right)^{\mu}.
\end{equation}
Define the mapping $F$
\begin{equation}\label{eq2.3}
    \begin{array}{rcl}
   F: \Omega^{B^{d_0}}(\mu) & \longrightarrow   & \Omega^{B^{d_0}}(\mu), \\
     (z,w)            & \longmapsto   & (\phi(z),\psi(z)w).
  \end{array}
\end{equation}
Then $F$ is an isometric automorphism of $(\Omega^{B^{d_0}}(\mu),g(\mu))$,
that is
\begin{equation}\label{eq2.4}
   \partial\overline{\partial}(\Phi(F(z,w)))=\partial\overline{\partial}(\Phi(z,w)).
\end{equation}
 }\end{Lemma}
\begin{proof}[Proof]
From \cite{RWYZ}, we know that $F$ is an automorphism of $\Omega^{B^{d_0}}(\mu)$, and
\begin{equation}\label{eq2.5}
    N(\phi(z),\overline{\phi(z)})^p=J\phi(z)N(z,\overline{z})^p\overline{J\phi(z)},
\end{equation}
where $J\phi(z)$ is the Jacobian of $\phi$.

By  \eqref{eq2.2} and \eqref{eq2.5}, we have
\begin{eqnarray*}
% \nonumber to remove numbering (before each equation)
   & & N(\phi(z),\overline{\phi(z)})^{\mu} -\|\psi(z)w\|^2\\
   &=&  N(\phi(z),\overline{\phi(z)})^{\mu}\left(1-\frac{\|w\|^2}{N(z,\overline{z})^{\mu}}\right) \\
   &=&|J\phi(Z)|^{\frac{2\mu}{p}}\left(N(z,\overline{z})^{\mu}-\|w\|^2\right),
\end{eqnarray*}
which implies \eqref{eq2.4}.
\end{proof}
\begin{Lemma}\label{Le:2.2}{
Let $\Omega$ be the Cartan domain with  the generic norm
$N(z,\overline{\xi})$, the dimension $d$ and the genus $p$. Then we
have
\begin{equation}\label{eq2.6}
   \det\left(\frac{\partial^2\Phi}{\partial z_i\partial
   \overline{z_j}}\right)_{i,j=1}^n(z,w)=\frac{\mu^{d}C_{\Omega}N(z,\overline{z})^{\mu(d+1)-p}}{(N(z,\overline{z})^{\mu}-\|w\|^2)^{n+1}},
\end{equation}
where the function $\Phi(z,w)=-\log(N(z,\overline{z})^{\mu}-\|w\|^2)$,
$n=d+d_0$,  $z=(z_1,z_2,\cdots,z_d)$,
$w=(w_1,w_2,\cdots,w_{d_0})=(z_{d+1},z_{d+2},\cdots,z_n)$, and
{$C_{\Omega}=\left.\det(-\frac{\partial^2\log
N(z,\overline{z})}{\partial
z_i\partial\overline{z_j}})\right|_{z=0}$}.
 }\end{Lemma}

\begin{proof}[Proof]
It is well known that
\begin{equation}\label{eq2.7}
    \frac{(\frac{\sqrt{-1}}{2\pi}\partial\overline{\partial}\Phi)^n}{n!}=\det\left(\frac{\partial^2\Phi}{\partial z_i\partial
   \overline{z_j}}\right)_{i,j=1}^n\frac{\omega_0^n}{n!},
\end{equation}
where $\omega_0=\frac{\sqrt{-1}}{2\pi}\sum_{j=1}^ndz_j\wedge
d\overline{z_j}$.

From \eqref{eq2.4} and \eqref{eq2.7}, we get
\begin{equation}\label{eq2.8}
\det\left(\frac{\partial^2\Phi(F)}{\partial z_i\partial
   \overline{z_j}}\right)_{i,j=1}^n=\det\left(\frac{\partial^2\Phi}{\partial z_i\partial
   \overline{z_j}}\right)_{i,j=1}^n.
\end{equation}
By the identity
\begin{equation}\label{eq2.9}
\left(\frac{\partial^2\Phi(F)}{\partial z_i\partial
   \overline{z_j}}\right)_{i,j=1}^n(z,w)=\left(\frac{\partial F_j}{\partial z_i}\right)_{i,j=1}^n\left(\frac{\partial^2\Phi}{\partial z_i\partial
   \overline{z_j}}\right)_{i,j=1}^n(F(z,w))\left(\frac{\partial \overline{F_i}}{\partial\overline{
   z_j}}\right)_{i,j=1}^n
\end{equation}
 and \eqref{eq2.8}, we deduce
\begin{equation}\label{eq2.10}
\det\left(\frac{\partial^2\Phi}{\partial z_i\partial
   \overline{z_j}}\right)_{i,j=1}^n(z,w)=|JF(z,w)|^2\det\left(\frac{\partial^2\Phi}{\partial z_i\partial
   \overline{z_j}}\right)_{i,j=1}^n(F(z,w)),
\end{equation}
where
\begin{equation*}\label{eq2.11}
   JF(z,w):=\det\left(\frac{\partial F_j}{\partial
   z_i}\right)_{i,j=1}^n,
\end{equation*}
and
\begin{equation*}
\left(\frac{\partial F_j}{\partial z_i}\right)_{i,j=1}^n:=\left(
                                                           \begin{array}{cccc}
                                                            \frac{\partial F_1}{\partial z_1 }  &  \frac{\partial F_2}{\partial z_1 } & \cdots & \frac{\partial F_n}{\partial z_1 }   \\
                                                            \frac{\partial F_1}{\partial z_2 }  &  \frac{\partial F_2}{\partial z_2 } & \cdots & \frac{\partial F_n}{\partial z_2 } \\
                                                             \vdots                             &  \vdots                             & \vdots & \vdots \\
                                                            \frac{\partial F_1}{\partial z_n }  &  \frac{\partial F_2}{\partial z_n } & \cdots & \frac{\partial F_n}{\partial z_n } \\
                                                           \end{array}
                                                         \right).
\end{equation*}
Let $(\widetilde{z_0},\widetilde{w_0})=F(z_0,w_0)$, $(z_0,w_0)\in
\Omega^{B^{d_0}}(\mu)$.  By \eqref{eq2.1} and \eqref{eq2.3}, then
$(\widetilde{z_0},\widetilde{w_0})=\left(0,\frac{w_0}{N(z_0,\overline{z_0})^{\frac{\mu}{2}}}\right)$
and
\begin{equation}\label{eq2.12}
 |JF(z_0,w_0)|^2=|J\phi(z_0)|^2|\psi(z_0)|^{2d_0}.
\end{equation}
Using $N(0,z)=1$, \eqref{eq2.1},  \eqref{eq2.5}, \eqref{eq2.12} and \eqref{eq2.10},
we have
\begin{equation}\label{eq2.13}
    |JF(z_0,w_0)|^2=\frac{1}{N(z_0,\overline{z_0})^{p+\mu d_0}},
\end{equation}
and
\begin{equation}\label{eq2.13}
\det\left(\frac{\partial^2\Phi}{\partial z_i\partial
   \overline{z_j}}\right)_{i,j=1}^n(z_0,w_0)=\frac{1}{N(z_0,\overline{z_0})^{p+\mu d_0}}\det\left(\frac{\partial^2\Phi}{\partial z_i\partial
   \overline{z_j}}\right)_{i,j=1}^n(0,\widetilde{w_0}).
\end{equation}

Now we calculate $\det\left(\frac{\partial^2\Phi}{\partial
z_i\partial\overline{z_j}}\right)_{i,j=1}^n(0,w)$.

From
\begin{equation}\label{eq2.14}
  \left\{\begin{array}{ll}
           N(0,0) & =1, \\
           \left.\frac{\partial \log N(z,\overline{z})}{\partial z_i}\right|_{z=0} &= \left.\frac{\partial \log N(z,\overline{z})}{\partial \overline{z_i}}\right|_{z=0}=0 \quad (1\leq i \leq d ),

         \end{array}
  \right.
\end{equation}
we obtain
\begin{equation}\label{eq2.17}
\left\{  \begin{array}{ll}
    \left.\frac{\partial N(z,\overline{z})^{\mu}}{\partial z_i}\right |_{z=0} & =\left.\frac{\partial }{\partial z_i}\exp\{\mu\log N(z,\overline{z})\}\right |_{z=0}=\left.\mu N(z,\overline{z})^{\mu} \frac{\partial \log N(z,\overline{z})}{\partial  z_i}\right |_{z=0}=0,  \\
     \left.\frac{\partial N(z,\overline{z})^{\mu}}{\partial\overline{ z_i}}\right |_{z=0} & =0, \\
     \left.\frac{\partial^2N(z,\overline{z})^{\mu}}{\partial z_i \partial \overline{z_j}}\right |_{z=0} & =\left.\left\{\mu N(z,\overline{z})^{\mu} \frac{\partial^2 \log N(z,\overline{z})}{\partial  z_i \partial
   \overline{z_j}}+\mu \frac{\partial \log N(z,\overline{z}) }{\partial z_i}\frac{\partial  N(z,\overline{z})^{\mu} }{\partial
   \overline{z_j}}\right\}\right |_{z=0}=-\mu c_{ij}.
  \end{array}\right.
\end{equation}
{where $c_{ij}=-\left.\frac{\partial^2 \log
N(z,\overline{z})}{\partial z_i\partial
\overline{z_j}}\right|_{z=0}$}.
Therefore
\begin{equation}\label{eq2.20}
  \left\{\begin{array}{ll}
           \frac{ \partial^2\Phi}{\partial z_i \partial   \overline{z_j}}(0,w) & =\left.\left\{\frac{1}{(N(z,\overline{z})^{\mu}-\|w\|^2)^2}\frac{\partial N(z,\overline{z})^{\mu}}{\partial z_i}\frac{\partial N(z,\overline{z})^{\mu}}{\partial   \overline{z_j}}-\frac{1}{N(z,\overline{z})^{\mu}-\|w\|^2}\frac{\partial^2 N(z,\overline{z})^{\mu}}{\partial z_i \partial   \overline{z_j}}\right\}\right|_{z=0} \\
            & = \frac{\mu c_{ij}}{1-\|w\|^2}\quad (1\leq i, j\leq d),\\
          \frac{ \partial^2\Phi}{\partial w_i \partial \overline{w_j}}(0,w)  &= \left.\left\{\frac{\delta_{ij}}{N(z,\overline{z})^{\mu}-\|w\|^2}+\frac{\overline{w_i}w_j}{(N(z,\overline{z})^{\mu}-\|w\|^2)^2}\right\}\right|_{z=0} \\
            &= \frac{\delta_{ij}}{1-\|w\|^2}+\frac{\overline{w_i}w_j}{(1-\|w\|^2)^2}\quad (1\leq i, j\leq d_0), \\
          \frac{ \partial^2\Phi}{\partial z_i \partial\overline{w_j}}(0,w)  & =\left.\frac{-w_j}{(N(z,\overline{z})^{\mu}-\|w\|^2)^2}\frac{\partial N(z,\overline{z})^{\mu}}{\partial z_i}\right|_{z=0}=0\quad (1\leq i\leq d, 1\leq j\leq d_0), \\
          \frac{ \partial^2\Phi}{\partial w_i \partial\overline{z_j}}(0,w)  & =\left.\frac{-\overline{w_i}}{(N(z,\overline{z})^{\mu}-\|w\|^2)^2}\frac{\partial
          N(z,\overline{z})^{\mu}}{\partial\overline{z_j}}\right|_{z=0}=0\quad(1\leq i\leq d_0, 1\leq j\leq d).
         \end{array}
  \right.
\end{equation}
 The above results can be rewritten as
 \begin{equation}\label{eq2.24}
\left(\frac{\partial^2\Phi}{\partial z_i\partial
   \overline{z_j}}\right)_{i,j=1}^n(0,w)=\left(
                                           \begin{array}{cc}
                                             \frac{\mu}{1-\|w\|^2}C_d & 0 \\
                                             0                        & \frac{1}{1-\|w\|^2}I_{d_0}+\frac{1}{(1-\|w\|^2)^2}w^{\dag}w \\
                                           \end{array}
                                         \right),
 \end{equation}
where $I_{d_0}$ denotes  the $d_0\times d_0$ diagonal matrix with
its diagonal elements 1,   $w^{\dag}$ is  the conjugate transpose of
the row vector $w=(w_1,w_2,\cdots,w_{d_0})$, and
{$C_d=(c_{ij})_{i,j=1}^d$}.

From \eqref{eq2.24}, we have
\begin{equation}\label{eq2.25}
    \det\left(\frac{\partial^2\Phi}{\partial z_i\partial
   \overline{z_j}}\right)_{i,j=1}^n(0,w)=\frac{\mu^d\det C_d }{(1-\|w\|^2)^{d+d_0+1}}.
\end{equation}
Finally, by \eqref{eq2.13} and \eqref{eq2.25}, we have
\eqref{eq2.6}.
\end{proof}

\begin{Theorem}\label{Th:a.2}{Let $\Omega^{B^{d_0}}(\mu)$  be a Cartan-Hartogs
domain with the K\"{a}hler metric $g(\mu)$ (see \eqref{eq1.3}). Let
$\alpha>\max\{d+d_0,\frac{p-1}{\mu}\}$. Then the Bergman kernel
$K_{\alpha}(z,w;\overline{z},\overline{w})$ of the Hilbert space
$$\mathcal{H}_{\alpha}=\left\{ f\in \textmd{Hol}(\Omega^{B^{d_0}}(\mu))~\left|~ \int_{\Omega^{B^{d_0}}(\mu)}\right.|f|^2\exp\{-\alpha \Phi\}
\frac{\omega(\mu)^{d+d_0}}{(d+d_0)!}<+\infty\right\}$$ can be
written as
\begin{eqnarray}
% \nonumber to remove numbering (before each equation)
\nonumber          & & K_{\alpha}(z,w;\overline{z},\overline{w})  \\
\label{eq2.26}
                   &=&\left.\frac{\pi^{d+d_0}\chi_2(\alpha-d-d_0-1)}{C_{\Omega}\mu^d\chi_1(0)\chi_2(0)V(\Omega)V(B^{d_0})}\left(\frac{1}{N(z,\overline{z})}\right)^{\mu\alpha}
                       \chi_1(\mu(\alpha+t\frac{d}{dt})-p)\frac{1}{\left(1-t\frac{\|w\|^2}{N(z,\overline{z})^{\mu}}\right)^{\alpha-d}}\right|_{t=1},
\end{eqnarray}
where {$C_{\Omega}=\left.\det(-\frac{\partial^2\log
N(z,\overline{z})}{\partial
z_i\partial\overline{z_j}})\right|_{z=0}$},  $\chi_1, \chi_2$ and
$V(\Omega), V(B^{d_0})$ are Hua polynomials (see \eqref{1.3}) and
the volumes with respect to the Euclidean measure of $\Omega$,
$B^{d_0}$, respectively.

}\end{Theorem}
\begin{proof}[Proof] By  \eqref{eq2.6}, the inner product on
$\mathcal{H}_{\alpha}$ is given by
$$(f,g)=\frac{\mu^dC_{\Omega}}{\pi^{d+d_0}}\int_{\Omega^{B^{d_0}}(\mu)}f(z,w)\overline{g(z,w)}N(z,\overline{z})^{\mu(\alpha-d_0)-p}\left(1-\frac{\|w\|^2}{N(z,\overline{z})^{\mu}}\right)^{\alpha-(d+d_0+1)}dm(z)dm(w),$$
where $dm$ denotes the Euclidean measure.

For convenience, we set $\Omega_1=\Omega$, $\Omega_2=B^{d_0}$. Let
$r_i, a_i, b_i, d_i, p_i, \chi_i,(s)_\mathbf{\lambda}^{(i)}$ and
$N_i$ be rank, characteristic multiplicities, dimension, genus, Hua
polynomial,  generalized Pochhammer symbol and generic norm of the
irreducible bounded symmetric domain $\Omega_i$, $1\leq i\leq 2$.

Let $\mathcal{G}_i$ stand for the identity connected components of
 groups of biholomorphic self-maps of $\Omega_i \subset \mathbb{C}^{d_i}$, and $\mathcal{K}_i$
for stabilizer of the origin in $\mathcal{G}_i$, respectively. For
any $k=(k_1,k_2)\in \mathcal{K}:= \mathcal{K}_1\times\mathcal{K}_2$,
we define the action
 $$ \pi(k)f(z,w)\equiv f\circ k(z,w):=f(k_1\circ z,k_2\circ w)$$
 of $\mathcal{K}$, then the space $\mathcal{P}$
of holomorphic polynomials on $\mathbb{C}^{d_1}\times
\mathbb{C}^{d_2}$ admits the Peter-Weyl decomposition
$$\mathcal{P}=\bigoplus_{{\ell(\lambda)\leq r_1\atop \ell(\nu)\leq r_2}}\mathcal{P}^{(1)}_{\lambda}\otimes \mathcal{P}^{(2)}_{\nu},$$
where  spaces $\mathcal{P}^{(i)}_{\lambda}$ are
$\mathcal{K}_i$-invariant and irreducible subspaces of spaces
 of holomorphic polynomials on $\mathbb{C}^{d_i}(1\leq i \leq 2)$.

Since $\mathcal{H}_{\alpha}$ is invariant under the action of
$\mathcal{K}_1\times \mathcal{K}_2$, namely, $\forall k\in
\mathcal{K}_1\times \mathcal{K}_2$, $(\pi(k)f,\pi(k)g)=(f,g)$,
 $\mathcal{H}_{\alpha}$ admits
an irreducible decomposition (see \cite{FT})
$$\mathcal{H}_{\alpha}=\widehat{\bigoplus_{{\ell(\lambda)\leq r_1\atop \ell(\nu)\leq    r_2}}}\mathcal{P}^{(1)}_{\lambda}\otimes \mathcal{P}^{(2)}_{\nu},$$
where $\widehat{\bigoplus}$ denotes the orthogonal direct sum.

For every partition $\lambda$ of length $\leq r_i$, let
$K^{(i)}_{\lambda}(z,\overline{u})$ be the Bergman kernel of
$\mathcal{P}^{(i)}_{\lambda}$ with respect to \eqref{1.4}. By
Schur's lemma, there exist positive constants $c_{\lambda\nu}$ such
that
$c_{\lambda\nu}K^{(1)}_{\lambda}(z,\overline{z})K^{(2)}_{\nu}(w,\overline{w})$
are the reproducing kernels of $\mathcal{P}^{(1)}_{\lambda}\otimes
\mathcal{P}^{(2)}_{\nu}$ with respect to the above inner product
$(\cdot,\cdot)$. According to the definition of the reproducing kernel,
we have
\begin{eqnarray*}
% \nonumber to remove numbering (before each equation)
   & &\frac{\mu^dC_{\Omega}}{\pi^{d+d_0}}\int_{\Omega^{B^{d_0}}(\mu)}c_{\lambda\nu}K^{(1)}_{\lambda}(z,\overline{z})K^{(2)}_{\nu}(w,\overline{w})N(z,\overline{z})^{\mu(\alpha-d_0)-p}\left(1-\frac{\|w\|^2}{N(z,\overline{z})^{\mu}}\right)^{\alpha-(d+d_0+1)}dm(z)dm(w)  \\
   &=&\dim\mathcal{P}^{(1)}_{\lambda}\dim\mathcal{P}^{(2)}_{\nu}.
\end{eqnarray*}

Therefore, the Bergman kernel of $\mathcal{H}_{\alpha}$ can be
written as
\begin{equation}\label{e4.4}
    K_{\alpha}(z,w;\overline{z},\overline{w})=\sum_{{\ell(\lambda)\leq r_1\atop \ell(\nu)\leq
    r_2}}\frac{\dim\mathcal{P}^{(1)}_{\lambda}\dim\mathcal{P}^{(2)}_{\nu}}{<K^{(1)}_{\lambda}(z,\overline{z})K^{(2)}_{\nu}(w,\overline{w})>}K^{(1)}_{\lambda}(z,\overline{z})K^{(2)}_{\nu}(w,\overline{w}),
\end{equation}
where $<f>$ denotes integral
$$\frac{\mu^dC_{\Omega}}{\pi^{d+d_0}}\int_{\Omega^{B^{d_0}}(\mu)}f(z,w)N(z,\overline{z})^{\mu(\alpha-d_0)-p}\left(1-\frac{\|w\|^2}{N(z,\overline{z})^{\mu}}\right)^{\alpha-(d+d_0+1)}dm(z)dm(w).$$

If $\mu\alpha-p>-1$ and $\alpha-d-d_0-1>-1$ (namely,
$\alpha>\max\{d+d_0,\frac{p-1}{\mu}\}$), combining (see \cite{F},
\cite{FS})
\begin{equation}\label{eq}
    \int_{\Omega}K_{\lambda}(z,\overline{z})N(z,\overline{z})^sdm(z)=\frac{\dim \mathcal{P}_{\lambda}}{(p+s)_{\lambda}}
    \int_{\Omega}N(z,\overline{z})^sdm(z)
\end{equation}
for $ s>-1 $  and \eqref{1.2},  we have
\begin{eqnarray}
% \nonumber to remove numbering (before each equation)
\nonumber   & & \frac{\mu^dC_{\Omega}}{\pi^{d+d_0}}\int_{\Omega^{B^{d_0}}(\mu)}K^{(1)}_{\lambda}(z,\overline{z})K^{(2)}_{\nu}(w,\overline{w})N(z,\overline{z})^{\mu(\alpha-d_0)-p}\left(1-\frac{\|w\|^2}{N(z,\overline{z})^{\mu}}\right)^{\alpha-(d+d_0+1)}dm(z)dm(w)  \\
\nonumber   &=&
\frac{\mu^dC_{\Omega}}{\pi^{d+d_0}}\int_{\Omega}K^{(1)}_{\lambda}(z,\overline{z})N(z,\overline{z})^{\mu(\alpha+\nu)-p}
dm(z)
\int_{B^{d_0}} K^{(2)}_{\nu}(w,\overline{w})(1-\|w\|^2)^{\alpha-(d+d_0+1)}dm(w)\\
\label{e4.5}
&=&\frac{\mu^{d}C_{\Omega}\chi_1(0)\chi_2(0)V(\Omega)V(B^{d_0})}{\pi^{d+d_0}\chi_1(\mu(\alpha+\nu)-p)\chi_2(\alpha-(d+d_0+1))}
\frac{\dim\mathcal{P}^{(1)}_{\lambda}\dim\mathcal{P}^{(2)}_{\nu}}{(\mu(\alpha+\nu))_{\lambda}^{(1)}(\alpha-d)_{\nu}^{(2)}}.
\end{eqnarray}

Combing  \eqref{e4.4}, \eqref{e4.5} and \eqref{1.7}£¬ we get
\begin{eqnarray}
% \nonumber to remove numbering (before each equation)
\nonumber   & & K_{\alpha}(z,w;\overline{z},\overline{w}) \\
\nonumber    &=&
  \sum_{{\ell(\lambda)\leq r_1\atop \ell(\nu)\leq r_2}}
  c\chi_1(\mu(\alpha+\nu)-p)(\mu(\alpha+\nu))_{\lambda}^{(1)}(\alpha-d)_{\nu}^{(2)}
    K^{(1)}_{\lambda}(z,\overline{z})K^{(2)}_{\nu}(w,\overline{w})   \\
\nonumber   &=&
c\sum_{{\ell(\nu)\leq r_2}} \chi_1(\mu(\alpha+\nu)-p)(\alpha-d)_{\nu}^{(2)} K^{(2)}_{\nu}(w,\overline{w})\frac{1}{N(z,\overline{z})^{\mu(\alpha+\nu)}}   \\
\nonumber   &=&
\frac{c}{{N(z,\overline{z})^{\mu\alpha}}}\sum_{{\ell(\nu)\leq
r_2}} \chi_1(\mu(\alpha+\nu)-p)(\alpha-d)_{\nu}^{(2)} K^{(2)}_{\nu}\left(\frac{w}{N(z,\overline{z})^{\mu}},\overline{w}\right) \\
\nonumber   &=&
\left.\frac{c}{{N(z,\overline{z})^{\mu\alpha}}}\sum_{{\ell(\nu)\leq
r_2}} \chi_1(\mu(\alpha+t\frac{d}{dt})-p)(\alpha-d)_{\nu}^{(2)} K^{(2)}_{\nu}\left(\frac{tw}{N(z,\overline{z})^{\mu}},\overline{w}\right)\right|_{t=1}   \\
\nonumber   &=&
\left.\frac{c}{{N(z,\overline{z})^{\mu\alpha}}}
\chi_1(\mu(\alpha+t\frac{d}{dt})-p)\sum_{{\ell(\nu)\leq
r_2}}(\alpha-d)_{\nu}^{(2)}K^{(2)}_{\nu}\left(\frac{tw}{N(z,\overline{z})^{\mu}},\overline{w}\right)\right|_{t=1}\\
\nonumber
&=&\left.\frac{c}{{N(z,\overline{z})^{\mu\alpha}}}
\chi_1(\mu(\alpha+t\frac{d}{dt})-p)\frac{1}{\left(1-\frac{t\|w\|^2}{N(z,\overline{z})^{\mu}}\right)^{\alpha-d}}\right|_{t=1},
\end{eqnarray}
where $$c=\frac{\pi^{d+d_0}\chi_2(\alpha-(d+d_0+1))}{\mu^{d}C_{\Omega}\chi_1(0)\chi_2(0)V(\Omega)V(B^{d_0})},$$
which completes the proof.
\end{proof}

In order to simplify  \eqref{eq2.26}, we need Lemma \ref{Le:2.1.1} below.

\begin{Lemma}\textup{(see \cite{F})}\label{Le:2.1.1}{
Let $\varphi(x)$ be a polynomial in $x$ of degree $n$ and let $Z$ be
a matrix of order $m$. Let $t$ be a real variable such that
$||tZ||<1$, {where $||Z||$ denotes the norm of  $Z$}. For a real
number $n_0$, take $x_0=-mn_0$. Then we have
\begin{equation}\label{2.4}
\varphi(t\frac{d}{dt})\frac{1}{{\det}(I-tZ)^{n_0}}=\frac{1}{{\det}(I-tZ)^{n_0}}\sum_{k=0}^n\frac{D^k\varphi(x_0)}{k!}\sum_{|\lambda|=k}\frac{|\lambda|!}{z_{\lambda}}n_0^{\ell(\lambda)}p_{\lambda}(\frac{1}{I-tZ}),
\end{equation}
where
$$\lambda=(1^{m_1(\lambda)} 2^{m_2(\lambda)}\cdots )\quad (m_i(\lambda)\geq 0),$$
$$ |\lambda|:=\sum_{i}im_i(\lambda),\quad
\ell(\lambda):=\sum_{i}m_i(\lambda),\quad
z_{\lambda}:=\prod_{i}i^{m_i(\lambda)}m_i(\lambda)!,$$
$$ p_{\lambda}(Z):=\prod_{i}(\textup{Tr}Z^i)^{m_i(\lambda)},\quad D^k\varphi(x_0)=\sum_{j=0}^k{k\choose
j}(-1)^j\varphi(x_0-j) .$$
 }\end{Lemma}

Combing Theorem \ref{Th:a.2} and Lemma \ref{Le:2.1.1}, we obtain the
explicit expression of the Bergman kernel $K_{\alpha}$ of the
Hilbert space $\mathcal{H}_{\alpha}$ as follows.

\begin{Theorem}\label{Th:2.2}
Assume that
\begin{equation}\label{eq2.29}
   \widetilde{\chi}(x):=\chi_1(\mu x-p)\equiv \prod_{j=1}^r\left(\mu
   x-p+1+(j-1)\frac{a}{2}\right)_{1+b+(r-j)a}.
\end{equation}
Let $D^k\widetilde{\chi}(x)$ be the  k-order difference of
$\widetilde{\chi}$ at $x$, that is
\begin{equation}\label{eq2.29.0}
D^k\widetilde{\chi}(x)=\sum_{j=0}^k{k\choose
j}(-1)^j\widetilde{\chi}(x-j).
\end{equation}
Then \eqref{eq2.26} can be rewritten as
\begin{eqnarray}
% \nonumber to remove numbering (before each equation)
\nonumber          & & K_{\alpha}(z,w;\overline{z},\overline{w})  \\
\label{eq2.26.1}
                   &=&\frac{\pi^{d+d_0}}{\mu^dC_{\Omega}\chi_1(0)\chi_2(0)V(\Omega)V(B^{d_0})}\left(\frac{1}{N(z,\overline{z})}\right)^{\mu\alpha}
                       \sum_{k=0}^d\frac{D^k\widetilde{\chi}(d)}{k!}\frac{(\alpha-d-d_0)_{k+d_0}}{\left(1-\frac{\|w\|^2}{N(z,\overline{z})^{\mu}}\right)^{\alpha-d+k}}.
\end{eqnarray}
\end{Theorem}
\begin{proof}[Proof]Let $x_0=d-\alpha$. By
$$\chi_1(\mu(\alpha+x)-p)|_{x=x_0-j}=\chi_1(\mu(d-j)-p)=\widetilde{\chi}(d-j),$$
 we have
 $$D^k(\chi_1(\mu(\alpha+x)-p))|_{x=x_0}=D^k\widetilde{\chi}(d).$$
 Using \eqref{2.4} and
 \begin{equation*}
    (x)_k=\sum_{|\lambda|=k}\frac{|\lambda|!}{z_{\lambda}}x^{\ell(\lambda)},
 \end{equation*}
 we have
\begin{equation}\label{eq2.31}
 \chi_1(\mu(\alpha+t\frac{d}{dt})-p)\frac{1}{(1-tz)^{\alpha-d}}=\frac{1}{(1-tz)^{\alpha-d}}\sum_{k=0}^d\frac{D^k\widetilde{\chi}(d)}{k!}\frac{(\alpha-d)_k}{(1-tz)^k}.
\end{equation}
For the Cartan domain $B^{d_0}$, its Hua polynomial
\begin{equation*}
    \chi_2(x)=(x+1)_{d_0}.
\end{equation*}
Then
\begin{equation}\label{eq2.32}
    \chi_2(\alpha-d-d_0-1)(\alpha-d)_k=(\alpha-d-d_0)_{d_0+k}.
\end{equation}
From \eqref{eq2.31} and \eqref{eq2.32}, we get \eqref{eq2.26.1}.
\end{proof}

\setcounter{equation}{0}
\section{The Rawnsley's $\varepsilon$-functions for $\Omega^{B^{d_0}}(\mu)$  with the canonical metric $g(\mu)$}

In this section we give the explicit expression of  the Rawnsley's
$\varepsilon$-function and the coefficients $a_1,a_2$ of its
expansion for the Cartan-Hartogs domain
$(\Omega^{B^{d_0}}(\mu),g(\mu))$.
\begin{Theorem}\label{Th:2.3}{
Let $\alpha>\max\{d+d_0,\frac{p-1}{\mu}\}$. Then the Rawnsley's
$\varepsilon$-function associated to
$(\Omega^{B^{d_0}}(\mu),g(\mu))$ can be written as
\begin{equation}\label{eq2.37}
    \varepsilon_{\alpha}(z,w)=\frac{1}{\mu^d}\sum_{k=0}^d\frac{D^k\widetilde{\chi}(d)}{k!}\left(1-\frac{\|w\|^2}
    {N(z,\overline{z})^{\mu}}\right)^{d-k}(\alpha-d-d_0)_{k+d_0}
\end{equation}
$($see \eqref{eq2.29} and \eqref{eq2.29.0} for the definition of the
functions $\widetilde{\chi}(x)$ and $D^k\widetilde{\chi}(x)$
respectively$)$. }\end{Theorem}
\begin{proof}[Proof]
By \eqref{eq2.26.1}, we have
\begin{eqnarray*}
% \nonumber to remove numbering (before each equation)
   & & \exp\{-\alpha    \Phi(z,w)\}K_{\alpha}(z,w;\overline{z},\overline{w}) \\
   &=&\frac{\pi^{d+d_0}}{\mu^dC_{\Omega}\chi_1(0)\chi_2(0)V(\Omega)V(B^{d_0})}\sum_{k=0}^d\frac{D^k\widetilde{\chi}(d)}{k!}\left(1-\frac{\|w\|^2}{N(z,\overline{z})^{\mu}}\right)^{d-k}(\alpha-d-d_0)_{k+d_0}.
\end{eqnarray*}
From \cite{Hua} and \cite{Ko}, we have
$$V(\Omega)=\frac{\pi^d}{C_{\Omega}\chi_1(0)}.$$
Since $C_{B^{d_0}}=1$, it follows that
$$V(B^{d_0})=\frac{\pi^{d_0}}{C_{B^{d_0}}\chi_2(0)}=\frac{\pi^{d_0}}{\chi_2(0)}.$$
Therefore, we obtain \eqref{eq2.37}.
\end{proof}
\begin{Corollary}\label{Co:2.1}{
 The coefficients $a_1$ and $a_2$ of the expansion of the Rawnsley's $\varepsilon$-function $\varepsilon_{\alpha}$, that
is, the coefficients of $\alpha^{d+d_0-1}$ and $ \alpha^{d+d_0-2}$
in
 \eqref{eq2.37} respectively,  are given by
\begin{equation}\label{eq2.38}
   a_1(z,w)=\frac{1}{\mu^d}\frac{D^{d-1}\widetilde{\chi}(d)}{(d-1)!}\left(1-\frac{\|w\|^2}{N(z,\overline{z})^{\mu}}\right)-\frac{(d+d_0)(d+d_0+1)}{2},
\end{equation}
\begin{eqnarray}
% \nonumber to remove numbering (before each equation)
\nonumber  a_2(z,w) &=& \frac{1}{\mu^d}\frac{D^{d-2}\widetilde{\chi}(d)}{(d-2)!}\left(1-\frac{\|w\|^2}{N(z,\overline{z})^{\mu}}\right)^2-  \\
\nonumber           & & \frac{1}{\mu^d}\frac{D^{d-1}\widetilde{\chi}(d)}{(d-1)!}\left(\frac{(d+d_0)(d+d_0+1)}{2}-1\right)
\left(1-\frac{\|w\|^2}{N(z,\overline{z})^{\mu}}\right)+ \\
\label{eq2.39} & &\frac{1}{24}(d+d_0-1)(d+d_0)(d+d_0+1)(3(d+d_0)+2)
\quad    (d\geq 2),
\end{eqnarray}
\begin{eqnarray}
% \nonumber to remove numbering (before each equation)
\nonumber   a_2(z,w) &=& -\frac{\mu-1}{\mu}\left(\frac{(1+d_0)(2+d_0)}{2}-1\right)\left(1-\frac{\|w\|^2}{N(z,\overline{z})^{\mu}}\right)+ \\
\label{eq2.39.1}     & & \frac{1}{24}d_0(d_0+1)(d_0+2)(3d_0+5)
 \quad (d=1).
\end{eqnarray}
 }\end{Corollary}
\begin{proof}[Proof]
Let
\begin{equation}\label{eq2.40}
    (\alpha-d-d_0)_{d_0+k}=\sum_{j=0}^{d_0+k}c_{d_0+k,j}\alpha^j.
\end{equation}
Substituting \eqref{eq2.40} into \eqref{eq2.37},  we obtain
\begin{equation}\label{eq2.41}
   \varepsilon_{\alpha}(z,w)=\sum_{j=0}^{d+d_0}\alpha^j\sum_{k=\max(j-d_0,0)}^d
  \frac{c_{d_0+k,j}}{\mu^d}
  \frac{D^k\widetilde{\chi}(d)}{k!}\left(1-\frac{\|w\|^2}{N(z,\overline{z})^{\mu}}\right)^{d-k},
\end{equation}
which implies
\begin{equation}\label{eq2.42}
a_j(z,w)=\sum_{k=\max(d-j,0)}^d
\frac{c_{d_0+k,d+d_0-j}}{\mu^d}\frac{D^k\widetilde{\chi}(d)}{k!}\left(1-\frac{\|w\|^2}{N(z,\overline{z})^{\mu}}\right)^{d-k}.
\end{equation}
By
\begin{eqnarray*}
% \nonumber to remove numbering (before each equation)
  (\alpha-d-d_0)_{d+d_0}   &=& \prod_{k=1}^{d+d_0}(\alpha-k), \\
  (\alpha-d-d_0)_{d+d_0-1} &=& \prod_{k=2}^{d+d_0}(\alpha-k), \\
  (\alpha-d-d_0)_{d+d_0-2} &=& \prod_{k=3}^{d+d_0}(\alpha-k),
\end{eqnarray*}
we have
\begin{equation}\label{eq2.44}
c_{d+d_0-1,d+d_0-1}=c_{d+d_0-2,d+d_0-2}=1,
\end{equation}
\begin{equation}\label{eq2.45}
c_{d+d_0,d+d_0-1}=-\sum_{k=1}^{d+d_0}k=-\frac{(d+d_0)(d+d_0+1)}{2},
\end{equation}
\begin{equation}\label{eq2.46}
c_{d+d_0-1,d+d_0-2}=-\sum_{k=2}^{d+d_0}k=-\frac{(d+d_0)(d+d_0+1)}{2}+1,
\end{equation}

\begin{eqnarray}
% \nonumber to remove numbering (before each equation)
\nonumber  c_{d+d_0,d+d_0-2} &=& \sum_{1\leq i<j\leq
                                 d+d_0}ij=\frac{1}{2}\left\{\left(\sum_{k=1}^{d+d_0}k\right)^2-\sum_{k=1}^{d+d_0}k^2\right\} \\
\label{eq2.47}&=&\frac{1}{24}(d+d_0-1)(d+d_0)(d+d_0+1)(3(d+d_0)+2).
\end{eqnarray}

For $d>1$, substituting \eqref{eq2.44}, \eqref{eq2.45},
\eqref{eq2.46} and \eqref{eq2.47} into \eqref{eq2.42}, we obtain
\eqref{eq2.38} and \eqref{eq2.39}.

For $d=1$, we have
\begin{equation}\label{eq2.48}
    \widetilde{\chi}(x)=\mu x-1.
\end{equation}
Thus, by \eqref{eq2.37}, we get \eqref{eq2.39.1}.
\end{proof}

In order to calculate $D^{d-1}\widetilde{\chi}$ and
$D^{d-2}\widetilde{\chi}$, we need Lemma \ref{Le:2.8} and
\ref{Le:2.9} below.

\begin{Lemma}\label{Le:2.8}{
Let $\widetilde{\chi}(x):= \prod_{j=1}^r\left(\mu
   x-p+1+(j-1)\frac{a}{2}\right)_{1+b+(r-j)a}=\sum_{j=0}^dc_jx^{d-j}$. Then
\begin{equation}\label{eq2.49}
    c_0=\mu^d,
\end{equation}
\begin{equation}\label{eq2.50}
   c_1=-\frac{1}{2}\mu^{d-1}dp,
\end{equation}

\begin{eqnarray}
% \nonumber to remove numbering (before each equation)
\nonumber c_2  &=& \frac{1}{2}\mu^{d-2}\left\{\frac{d^2p^2}{4}-\frac{r(p-1)p(2p-1)}{6}+\frac{r(r-1)a(3p^2-3p+1)}{12}-\right. \\
\label{eq2.51}&&\left.\frac{(r-1)r(2r-1)a^2(p-1)}{24}+\frac{r^2(r-1)^2a^3}{48}\right\}.
\end{eqnarray}
 }\end{Lemma}

 \begin{proof}[Proof]
 The coefficients of $x^{d}$, $ x^{d-1}$  and $ x^{d-2}$ in the function $\widetilde{\chi}(x)$ (see \eqref{eq2.29}),  respectively, are
given by
\begin{eqnarray}
% \nonumber to remove numbering (before each equation)
\nonumber       c_0 &=& \mu^d, \\
\label{eq2.52}  c_1 &=& -\mu^d\sum_{j=1}^r\sum_{i=1}^{1+b+(r-j)a}\frac{1}{\mu}\left(p-i-(j-1)\frac{a}{2}\right), \\
\label{eq2.53}  c_2 &=&\frac{1}{2}\mu^d\left\{\left(-
\frac{c_1}{\mu^d}
  \right)^2-\sum_{j=1}^r\sum_{i=1}^{1+b+(r-j)a}\left(\frac{1}{\mu}\left(p-i-(j-1)\frac{a}{2}\right)\right)^2\right\}.
\end{eqnarray}
Using \eqref{1.1} and formulas
$$ \sum_{k=1}^nk^2=\frac{n(n+1)(2n+1)}{6},~~ \sum_{k=1}^nk^3=\frac{n^2(n+1)^2}{4},$$
we have
\begin{eqnarray}
% \nonumber to remove numbering (before each equation)
\nonumber   \sum_{j=1}^r\sum_{i=1}^{1+b+(r-j)a}\left(p-i-(j-1)\frac{a}{2}\right) &=&\sum_{j=1}^r\sum_{k=1+(j-1)\frac{a}{2}}^{p-1-(j-1)\frac{a}{2}}k  \\
\nonumber   &=& \sum_{j=1}^r\frac{1}{2}p(p-1-(j-1)a)\\
\nonumber   &=& \frac{1}{2}p\sum_{j=1}^r(b+1+(r-j)a) \\
\nonumber   &=& \frac{1}{2}p\left((1+b)r+\frac{r(r-1)}{2}a\right)\\
\label{eq2.54} &=&\frac{1}{2}pd,
\end{eqnarray}
\begin{eqnarray}
% \nonumber to remove numbering (before each equation)
\nonumber   & &  \sum_{j=1}^r\sum_{i=1}^{1+b+(r-j)a}\left(p-i-(j-1)\frac{a}{2}\right)^2 \\
\nonumber   &=&  \sum_{j=1}^r\sum_{k=1+(j-1)\frac{a}{2}}^{p-1-(j-1)\frac{a}{2}}k^2\\
\nonumber   &=& \frac{1}{24}\sum_{j=0}^{r-1}\left\{(2p-2-ja)(2p-ja)(2p-1-ja)-ja(2+ja)(1+ja)\right\} \\
\nonumber   &=& \frac{1}{24}\sum_{j=0}^{r-1}\left\{(2p-2)(2p-1)2p-(12p^2-12p+4)aj+6(p-1)a^2j^2-2a^3j^3\right\} \\
\nonumber   &=& \frac{r(p-1)p(2p-1)}{6}-\frac{r(r-1)a(3p^2-3p+1)}{12}+\\
\label{eq2.55}&&\frac{(r-1)r(2r-1)a^2(p-1)}{24}-\frac{r^2(r-1)^2a^3}{48}.
\end{eqnarray}
Combining \eqref{eq2.52}, \eqref{eq2.53}, \eqref{eq2.54} and
\eqref{eq2.55}, we get \eqref{eq2.50} and \eqref{eq2.51}.
\end{proof}

\begin{Lemma}\label{Le:2.9}{
For any  polynomial $f(x)$ in real variable $x$, take
$Df(x):=f(x)-f(x-1)$. Let $A_d=D^{d-1}x^d$, $B_d=D^{d-2}x^d$. Then
we have
\begin{eqnarray}
% \nonumber to remove numbering (before each equation)
\label{eq2.56} A_d  &=& \frac{d!}{2}(2x-d+1)\quad  (d\geq 1),  \\
\label{eq2.57} B_d  &=&
\frac{d!}{24}\left\{12x^2-12(d-2)x+3d^2-11d+10\right\}  \quad (d\geq
2).
\end{eqnarray}
 }\end{Lemma}

 \begin{proof}[Proof]
Firstly, we have the recurrence relations£º
\begin{eqnarray}
% \nonumber to remove numbering (before each equation)
\nonumber A_d  &=& D^{d-2}(Dx^d)  \\
\nonumber      &=& D^{d-2}\left(\sum_{j=0}^{d-1}(-1)^{d+1-j}{d\choose j}x^j\right) \\
\nonumber      &=& d D^{d-2}x^{d-1}-\frac{d(d-1)}{2} D^{d-2}x^{d-2} \\
\label{eq2.58} &=& dA_{d-1}-\frac{d!}{2}
\end{eqnarray}
for $d>1$ and
\begin{eqnarray}
% \nonumber to remove numbering (before each equation)
\nonumber B_d &=& D^{d-3}(Dx^d) \\
\nonumber     &=& D^{d-3}\left(\sum_{j=0}^{d-1}(-1)^{d+1-j}{d\choose j}x^j\right)\\
\nonumber     &=& d D^{d-3}x^{d-1}-\frac{d(d-1)}{2} D^{d-3}x^{d-2} +\frac{d(d-1)(d-2)}{6}D^{d-3}x^{d-3}\\
\label{eq2.59}&=& dB_{d-1}-\frac{d(d-1)}{2}A_{d-2}+\frac{d!}{6}
\end{eqnarray}
for $d>2$.  Now,  by solving difference equation
\begin{equation}\label{eq2.60}
    \left\{\begin{array}{cl}
             A_d= & dA_{d-1}-\frac{d!}{2}, \\
             B_d= & dB_{d-1}-\frac{d(d-1)}{2}A_{d-2}+\frac{d!}{6},\\
             A_1= & x, \\
             B_2= & x^2,
           \end{array}
    \right.
\end{equation}
we obtain \eqref{eq2.56} and \eqref{eq2.57}.
\end{proof}

Lemma \ref{Le:2.8} and Lemma \ref{Le:2.9} imply the following results.
\begin{Lemma}\label{Le:2.10}{
Suppose that $D^{d-1}\widetilde{\chi}(d)$ and
$D^{d-2}\widetilde{\chi}(d)$ are defined by \eqref{eq2.29} and \eqref{eq2.29.0}. Then
we have
\begin{equation}\label{eq2.61}
    \frac{D^{d-1}\widetilde{\chi}(d)}{(d-1)!}=\frac{d\mu^{d-1}}{2}(\mu(d+1)-p) \quad  (d\geq 1),
\end{equation}
\begin{equation}\label{eq2.62}
    \frac{D^{d-2}\widetilde{\chi}(d)}{(d-2)!}=\mu^{d-2}\left\{\frac{1}{24}(d-1)d(d+1)(3d+10)\mu^2-\frac{1}{4}p(d-1)d(d+2)\mu+\frac{1}{2}
    \widetilde{c_2}\right\} \quad
    (d\geq 2),
\end{equation}
where
\begin{eqnarray}
% \nonumber to remove numbering (before each equation)
\nonumber \widetilde{c_2} &=& \frac{d^2p^2}{4}-\frac{r(p-1)p(2p-1)}{6}+\frac{r(r-1)a(3p^2-3p+1)}{12}- \\
 \label{eq2.63}           & & \frac{(r-1)r(2r-1)a^2(p-1)}{24}+\frac{r^2(r-1)^2a^3}{48}.
\end{eqnarray}
 }\end{Lemma}
\begin{proof}[Proof]From
$\widetilde{\chi}(x)=c_0x^d+c_1x^{d-1}+c_2x^{d-2}+\cdots+c_n$, we
get
\begin{equation}\label{eq2.64}
    \left\{\begin{array}{cl}
          D^{d-1}\widetilde{\chi}(x) =   & c_0A_d(x)+c_1(d-1)!,\\
          D^{d-2}\widetilde{\chi}(x) =   &          c_0B_d(x)+c_1A_{d-1}(x)+c_2(d-2)!.
           \end{array}
    \right.
\end{equation}
Let $x=d$, substituting \eqref{eq2.49}, \eqref{eq2.50},
\eqref{eq2.51}, \eqref{eq2.56} and \eqref{eq2.57} into
\eqref{eq2.64}, we obtain \eqref{eq2.61} and \eqref{eq2.62}.
\end{proof}

 \setcounter{equation}{0}
\section{The proof of Theorem \ref{Th:a.1.2}}
\begin{proof}[The proof of Theorem \ref{Th:a.1.2}]

If the dimension of $\Omega$ is 1 (i.e., $d=1$), from
\eqref{eq2.39.1}, we get that $a_2$ is constant if and only if
$\mu=1$, that is, $
 (\Omega^{B^{d_0}}(\mu), g(\mu)) $ is biholomorphically isometric to the complex
hyperbolic space $({B}^{1+d_0}, g_{hyp})$.

If the dimension of $\Omega$ is larger than 1 (i.e., $d>1$), it
follows from \eqref{eq2.39} that the coefficient $a_2$ of the
expansion of the function $\varepsilon_{\alpha}$ associated to
$(\Omega^{B^{d_0}}(\mu), g(\mu)) $ is constant if and only if
\begin{equation}\label{eq3.1}
    \frac{D^{d-1}\widetilde{\chi}(d)}{(d-1)!}=
    \frac{D^{d-2}\widetilde{\chi}(d)}{(d-2)!}=0.
\end{equation}
From \eqref{eq2.61},  \eqref{eq2.62} and \eqref{eq2.63}, we get that
$a_2$ is constant if and only if
\begin{equation}\label{eq3.2}
   \mu=\frac{p}{d+1},
\end{equation}
and
\begin{eqnarray}
% \nonumber to remove numbering (before each equation)
\nonumber       & & 12(d+1)\left\{\frac{d^2p^2}{4}-\frac{r(p-1)p(2p-1)}{6}+\frac{r(r-1)a(3p^2-3p+1)}{12}\right.\\
\label{eq3.3}
&&\left.-\frac{(r-1)r(2r-1)a^2(p-1)}{24}+\frac{r^2(r-1)^2a^3}{48}
\right\}-(d-1)d(3d+2)p^2=0.
\end{eqnarray}

  (1) For the bounded symmetric domain
$\Omega_I(m,n)\; (1\leq m\leq n)$, its  rank $r=m$, the
characteristic multiplicities $a=2, b=n-m$, the dimension $d=mn$,
the genus $p=m+n$. By \eqref{eq3.3}, we obtain
\begin{equation}\label{eq3.4}
  2mn(m^2-1)(n^2-1)=0,
\end{equation}
that is, $r=m=1$.

(2) For the bounded symmetric domain  $\Omega_{II}(2n)\; (n\geq 2) $,
its rank $r=n$, the characteristic multiplicities $a=4, b=0$, the
dimension $d=n(2n-1)$, the genus $p=2(2n-1)$. By \eqref{eq3.3}, we
obtain
\begin{equation}\label{eq3.5}
    4n^2(8n^4 - 20n^3 + 10n^2 + 5n - 3)=0,
\end{equation}
which is not satisfied by any  positive integer $n$ with $n\geq 2. $

 (3) For the bounded symmetric domain
$\Omega_{II}(2n+1)\;(n\geq 2) $, its  rank $r=n$, the characteristic
multiplicities $a=4, b=2$, the dimension $d=n(2n+1)$, the genus
$p=4n$. By \eqref{eq3.3}, we obtain
\begin{equation}\label{eq3.6}
    4n(2n - 1)(2n + 1)^2(n - 1)(n + 1)=0.
\end{equation}
This equation has no positive integer solution $n$ with $n\geq 2.$

 (4) For the bounded symmetric domain $\Omega_{III}(n)\;(n\geq 2) $, its rank
$r=n$, the characteristic multiplicities $a=1, b=0$, the dimension
$d=n(n+1)/2$, the genus $p=n+1$. By \eqref{eq3.3}, we obtain
\begin{equation}\label{eq3.7}
   \frac{1}{8}n^2(n^4 + 5n^3 + 5n^2 - 5n - 6)=0.
\end{equation}
The equation has no positive integer solution $n$ with $n\geq 2.$

(5) For the bounded symmetric domain  $\Omega_{IV}(n)\;(n\geq 5) $, its
rank $r=2$, the characteristic multiplicities $a=n-2, b=0$, the
dimension $d=n$, the genus $p=n$. By \eqref{eq3.3}, we obtain
\begin{equation}\label{eq3.8}
   n(n^2 + n - 2)=0
\end{equation}
for $n\geq 5$, which is impossible.

(6) For the bounded symmetric domain  $\Omega_{\textrm{V}}(16) $,
its rank $r=2$, the characteristic multiplicities $a=6, b=4$, the
dimension $d=16$, the genus $p=12$. It means that \eqref{eq3.3} does
not hold in this case.

 (7) For the bounded symmetric domain  $\Omega_{\textrm{VI}}(27) $, its  rank
$r=3$, the characteristic multiplicities $a=8, b=0$, the dimension
$d=27$, the genus $p=18$. It is obvious that \eqref{eq3.3} does not
hold in this case.

Combing the above results, we get that $a_2$ is constant if and only
if the rank $r=1$ of the bounded symmetric domain $\Omega$ and
$\mu=\frac{p}{d+1}=1$, which implies that the Cartan-Hartogs domain
$(\Omega^{B^{d_0}}(\mu), g(\mu)) $ is
 biholomorphically isometric to the complex hyperbolic space $({B}^{d+d_0}, g_{hyp})$.
\end{proof}

\noindent\textbf{Acknowledgments}\quad The part of the work was
completed when the first author visited School of Mathematics and
Statistics at Wuhan University during 2013, and he wishes to thank
the School for its kind hospitality. In addition, the authors would
like to thank the referees for many helpful suggestions. The first
author was supported by the Scientific Research Fund of Sichuan
Provincial Education Department (No.11ZA156), and the second author
was supported by the National Natural Science Foundation of China
(No.11271291).

%%%%%%%%%%%%%%%%%%%%%%%%%%%%%%%%%%%%%%%%%%%%%%%%%%%%%%%%%%%%%%%%

%%%%%%%%%%%%%%%%%%%%%%%%%%%%%%%%%%%%%%%%%%%%%%%%%%%%%%%%%%%%%%%%
\addcontentsline{toc}{section}{References}
\phantomsection
\renewcommand\refname{References}
\small{
}
%%%%%%%%%%%%%%%%%%%%%%%%%%%%%%%%%%%%%%%%%%%%%%%%%%%%%%%%%%%%%%%%
\clearpage
\end{document}